\documentclass[11pt,twoside]{amsart}

\usepackage{amsfonts}

\newtheorem{Theo1}{{Theorem}}
\newtheorem{Lemma1}{{Lemma}}[section]
\newtheorem{Def1}{{Definition}}
\newtheorem{Prop1}[Lemma1]{{Proposition}}
\newtheorem{Claim1}[Lemma1]{{Claim}}
\newtheorem{Rem1}[Lemma1]{{Remark}}
\newtheorem{Cor1}[Lemma1]{{Corollary}}
\newtheorem{Ex1}[Lemma1]{{Example}}

\newenvironment{Lemma}{\begin{Lemma1}}{\end{Lemma1}}
\newenvironment{Def}{\begin{Def1}\em}{\end{Def1}}
\newenvironment{Prop}{\begin{Prop1}}{\end{Prop1}}
\newenvironment{Rem}{\begin{Rem1}\rm}{\end{Rem1}}
\newenvironment{Theorem}{\begin{Theo1}}{\end{Theo1}}

\newenvironment{Cor}{\begin{Cor1}}{\end{Cor1}}

\newcommand{\uar}{\uparrow}
\newcommand{\dar}{\downarrow}
\newcommand{\lra}{\longrightarrow}
\newcommand{\lla}{\longleftarrow}
\newcommand{\ra}{\rightarrow}
\newcommand{\epi}{\lra \kern-.8em\ra}
\newcommand{\sdp}{\times\kern-.2em\vrule height1.1ex depth-.05ex}

\newcommand{\N}{{\mathbb N}}

\newcommand{\T}{{\mathbb T}}

\newcommand{\dickebox}{{\vrule height5pt width5pt depth0pt}}
\newcommand{\BA}{{\mathbf B}A}
\newcommand{\BB}{{\mathbf B}B}
\newcommand{\B}{{\mathbf B}}

 \normalbaselines
\title[Hochschild invariants of derived categories]{Fine
Hochschild invariants of derived categories\\ for symmetric
algebras}

\author{Alexander Zimmermann}
\address{Universit\'e de Picardie,\newline
Facult\'e de Math\'ematiques et CNRS UMR 6140, \newline 33 rue St Leu,
\newline
F-80039 Amiens Cedex 1, \newline France}
\email{alexander.zimmermann@u-picardie.fr}


\date{December 16, 2005; revised August 24, 2006}
\begin{document}

\begin{abstract}Let $A$ be a symmetric $k$-algebra over a
perfect field $k$. K\"ulshammer defined for any integer $n$
a mapping $\zeta_n$ on the degree $0$ Hochschild cohomology
and a mapping $\kappa_n$ on the degree $0$ Hochschild homology of $A$
as adjoint mappings of the respective $p$-power mappings with respect
to the symmetrizing bilinear form. In an earlier paper it is shown that
$\zeta_n$ is invariant under derived equivalences. In the present paper
we generalize the definition of $\kappa_n$ to higher Hochschild homology
and show the invariance of $\kappa$ and its generalization
under derived equivalences. This provides fine
invariants of derived categories.
\end{abstract}

\maketitle

\section*{Introduction}

Let $k$ be a commutative ring and let $A$ be a $k$-algebra which
is projective as a $k$-module. If $B$ is a second $k$-algebra and if
the derived categories $D^b(A)$ of bounded complexes of $A$-modules
and $D^b(B)$ of $B$-modules are equivalent as triangulated categories,
then  the Hochschild cohomology of $A$ is isomorphic to the
Hochschild cohomology of $B$
(Rickard \cite{Ri3}). Analogous statements hold for
the cyclic homology (Keller \cite{Kellercyclic}),
the $K$-theory (Thomason-Trobaugh \cite{thomason}),
the fact of being symmetric algebra (cf Rickard \cite{Ri3} for fields and
\cite{rogquest} in a more general situation) and others.
This is one of the reasons why the derived category $D^b(A)$
is now one of the main tools in representation theory.
Nevertheless, most of the invariants are quite difficult to compute, except
some small cases like the centre, or the Grothendieck group.
Therefore, it is usually quite hard to distinguish two derived categories.

A symmetric $k$-algebra $A$ is equipped with a symmetric non degenerate
bilinear form $(\;,\;):A\times A\lra k$. Denote by $KA$ the commutator
subspace, that is the $k$-linear space generated by the set of $ab-ba$
for $a,b\in A$. If $k$ is a perfect field
of characteristic $p>0$, then K\"ulshammer defined in \cite{Kueldeutsch}
$T_n(A)^\perp$ to be the orthogonal space to
the set of $x\in A$ so that $x^{p^n}$ falls into $KA$.
It turned out that $T_n(A)^\perp$ is a decreasing
sequence of ideals in the centre of $A$. If $A$ and $B$ are
symmetric $k$-algebras with equivalent derived categories, then the centres
of $A$ and $B$ are isomorphic, and
in \cite{kuelsquest} we showed that this isomorphism maps $T_n(A)^\perp$
to $T_n(B)^\perp$.
This fine ideal structure of the centre of the algebra $A$ gives valuable
and computable derived invariants of $A$.
In joint work with Thorsten Holm \cite{hz-tame} we are able to apply the
invariance of the ideals $T_n(A)^\perp$ to
tame blocks of group rings solving delicate questions whether
certain parameters in the defining relations of particular algebras
lead to different derived categories. Thorsten Holm and Andrzej
Skowro\'nski use this new fine invariant to classify all tame domestic
symmetric algebras up to derived equivalence \cite{Holmskow}.

Alejandro Adem asked during the 2005 Oberwolfach conference ''Cohomology of
finite groups: Interactions and applications'' if it is possible to
generalize K\"ulshammer's ideals of the centre of $A$ to
a derived invariant of higher degree Hochschild cohomology. The purpose of
this paper is to answer to this question.

K\"ulshammer shows in \cite{Kueldeutsch} that $T_n(A)^\perp$ is the image
of a certain mapping $\zeta_n:Z(A)\lra Z(A)$ which is defined by
$(z,a^{p^n})=(\zeta_n(z),a)^{p^n}$ for all $z\in Z(A)$ and all $a\in A$.
This is the way we view
$T_n(A)^\perp$ in \cite{kuelsquest}. Since $KA=Z(A)^\perp$ the dual operation
defines a mapping $\kappa_n:A/KA\lra A/KA$ by the equation
$(z^{p^n},a)=(z,\kappa_n(a))^{p^n}$ for all $z\in Z(A)$ and all $a\in A/KA$.
As a first result we show that $\kappa_n$ as well is a derived invariant,
observing that $HH_0(A)=A/KA$ which is known to be a derived invariant.

As a first step
we are going to show that for any symmetric algebra $A$ one gets a
non degenerate pairing
$$(\;,\;)_m:HH^m(A,A)\times HH_m(A,A)\lra k\;.$$
Since $HH_*(A,A)$ does not have a multiplicative structure, it seems to be
impossible to write down the defining relation for an analogue for
$\zeta_n^A$ on higher Hochschild cohomology. Nevertheless,
Hochschild cohomology is a graded commutative ring.
Suppose either $m$ is even and $p$ odd or $p$ even and $m$ arbitrary.
We are able to show for these parameters that the $p^n$-power mapping
$HH^{m}(A,A)\lra HH^{p^nm}(A,A)$ has a right adjoint
$\kappa_n^{(m),A}:HH_{p^nm}(A,A)\lra HH_{m}(A,A)$ with respect
to $(\;,\;)_m$ and $(\;,\;)_{p^nm}$ and moreover $\kappa_n^{(0)}=\kappa_n$.
As a main result we show that any derived equivalence $F$ of standard type
between $A$ and $B$ induces an isomorphism $HH_m(F)$ on the Hochschild
homology and this in turn conjugates $\kappa_n^{(m),A}$ to $\kappa_n^{(m),B}$.

Using a suggestion of Bernhard Keller we also study the $p$-power map
by the Gerstenhaber Lie structure of the Hochschild cohomology.
We show that this again is invariant under derived equivalences.
Nevertheless, there is no obvious reason why this $p$-power map should
be semilinear. Moreover, it is
only defined from degree $1$ Hochschild cohomology onwards. Furthermore,
we expect that this will be somewhat harder to compute in examples.
On the other hand this gives for $p=2$ a richer structure to the set
of derivations on $A$, and this structure is then a computable
derived invariant.

Our paper is organised as follows. In Section~\ref{Hochschildandderived}
we recall the basic constructions concerning
Hochschild (co-)homology and their invariance under derived
equivalences. In Section~\ref{kappa} we show the derived invariance of
K\"ulshammer's mapping $\kappa_n$ as well as some consequences for
derived equivalences between blocks of group rings. Section~\ref{cupkappa}
is devoted to the definition of the generalization of $\kappa_n$ to
higher Hochschild homology, and to show its invariance under derived
equivalence. In Section~\ref{Gerstenhaber} we first recall the Gerstenhaber
construction of a Lie algebra structure on the Hochschild cohomology, define
for all primes $p$ a structure of a restricted Lie algebra
on odd degree Hochschild homology and on the entire Hochschild
cohomology for $p=2$, and show its invariance under derived equivalence.

\section{Derived equivalences and Hochschild constructions revisited}
\label{Hochschildandderived}

\subsection{Basic constructions and definitions for Hochschild homology
and cohomology}\label{Hochschilddefinitions}

In order to fix notation and also for convenience of the reader we recall
in this section the basic notions for Hochschild cohomology
and homology. Most of the material in this section
can be found in Loday \cite{Loday}, Keller \cite{KellerGersten} and
Stasheff \cite{Stasheff}.

Let $k$ be a field and let $A$ be a $k$-algebra.
Consider a complex $\B A$ whose degree $n$ homogeneous component
is $(\B A)_n:=A^{\otimes (n+2)}$
and whose  differential
$d_n:(\B A)_{n}\lra (\B A)_{n-1}$ is
$$d_n(x_0\otimes\dots\otimes x_{n+1}):=
\sum_{j=0}^{n}(-1)^j
(x_0\otimes\dots\otimes x_{j-1}\otimes x_jx_{j+1}
\otimes x_{j+2}\otimes\dots\otimes x_{n+1}).$$
Then the complex $(\BA{},d)$ is a projective resolution of
$A$ as $A\otimes_kA^{op}$-module.

For a later application in Section~\ref{Gerstenhaber}
we need to extend this construction to
$$\B (A):= k\oplus A[1]\oplus (A\otimes A)[2]\oplus
(A\otimes A\otimes A)[3]\oplus\dots$$
where the brackets indicate the degrees of the components.
The mappings
\begin{eqnarray*}
\Delta_n:\B (A)_n&\lra &\oplus_{1\leq j\leq n}\B (A)_j\otimes \B(A)_{n-j}\\
x_1\otimes\dots\otimes x_{n}&\mapsto&
\sum_{j=0}^{n}(x_1\otimes\dots\otimes x_j)\otimes
(x_{j+1}\otimes\dots\otimes x_{n})
\end{eqnarray*}
(where it is understood that the boundary cases $j=0$ and $j=n$ correspond to
the element $1_k$ in the first or the last bracket)
compose to a graded co-algebra map $\Delta:\B(A){}\lra \B(A){}\otimes \B(A){}$
and $\B(A){}$ becomes a differential graded co-associative co-algebra.
Observe the shift by degree $2$ between $\B A$ and $\B(A)$.

We abbreviate in the sequel
$A^e:=A\otimes A^{op}$ and $A^*=Hom_k(A,k)$.

By definition, for all $A^e$-modules $M$ one puts
$$HH_n(A,M)= H_n(\B A\otimes_{A\otimes A^{op}}M)$$
and
$$HH^n(A,M)= H^n(Hom_{A\otimes A^{op}}(\B A,M)).$$

\subsection{Some facts on derived equivalences}

Suppose two finite dimensional $k$-algebras
$A$ and $B$ have equivalent derived categories
$D^b(A)\simeq D^b(B)$ as triangulated categories. Then, there is a
complex $Y$ in $D^b(B\otimes A^{op})$ and a complex $X$ in
$D^b(A\otimes B^{op})$ so that
$$F_Y:=Y\otimes_A^{\mathbb L}-:D^b(A)\lra D^b(B)$$
is an equivalence and so that
$$X\otimes_B^{\mathbb L}-:D^b(B)\lra D^b(A)$$
is an equivalence quasi-inverse to $F_Y$.
It is known that one may choose $X$ and $Y$ so that both
complexes are formed by projective modules if restricted to either side,
and that then the left derived tensor product can be replaced by the
ordinary tensor product.

Moreover, doing so,
$$F_Y^e:=Y\otimes_A^{\mathbb L}-\otimes_A^{\mathbb L}X:
D^b(A\otimes A^{op})\lra D^b(B\otimes B^{op})$$
is an equivalence of triangulated categories satisfying
$$F_Y^e(A)\simeq B\;.$$
In \cite{rogquest} we have shown that for quite general algebras,
in particular for finite dimensional algebras over fields $k$,
we have $$F_Y^e(Hom_k(A,k))\simeq Hom_k(B,k)\;.$$
Moreover, in \cite{kuelsquest} we have shown that for $k$ a field of
characteristic $p>0$ denoting by $k^{(n)}$
the $k$-vector space $k$ twisted by the $n$-th power of
the Frobenius automorphism $Fr$,
$$F_Y^e(Hom_k(A,k^{(n)}))\simeq Hom_k(B,k^{(n)}).$$

In particular, if $A\simeq Hom_k(A,k)$ in $A^e-mod$, then
$B\simeq Hom_k(B,k)$ in $B^e-mod$. Such algebras are called symmetric.

Another consequence (cf Rickard\cite{Ri3} or \cite{derbuch})
is that $F_Y^e$ (and therefore $F_Y$) induces an isomorphism
\begin{eqnarray*}
HH^m(A,A)=Ext^m_{A^e}(A,A)&=&Hom_{D^b(A^e)}(A,A[m])\\
&\stackrel{F^e}{\lra}&
Hom_{D^b(B^e)}(B,B[m])=HH^m(B,B)\;.
\end{eqnarray*}

Hochschild homology as well is an invariant of the derived category.
This seems to be well known, but I could not find a reference in the
literature, and therefore I include a proof. In particular, it is
important to have an explicit isomorphism which we will need for our
proof.
In fact, $\BA{}$ is a free resolution of $A$ in $A^e-mod$. Since $F^e_Y$
is an equivalence,
$$F_Y^e(\BA{})=Y\otimes_A\BA{}\otimes_AX$$
is a projective resolution of $B$ in $B^e-mod$. The mapping
\begin{eqnarray*}
\BA{}\otimes_{A^e}(X\otimes_BY)&\lra&
(Y\otimes_A\BA{}\otimes_AY)\otimes_{B^e}B\\
u\otimes (x\otimes y)&\mapsto&(y\otimes u\otimes x)\otimes 1
\end{eqnarray*}
is well defined as is easily seen.
Moreover, an inverse is
\begin{eqnarray*}
(Y\otimes_A\BA{}\otimes_AY)\otimes_{B^e}B&\lra&\BA{}\otimes_{A^e}(X\otimes_BY)\\
(y\otimes u\otimes x)\otimes b&\mapsto&u\otimes (xb\otimes y)=
                                    u\otimes (x\otimes by)
\end{eqnarray*}
We take the degree $m$ homology of the various complexes and
observe since $X\otimes_BY\simeq A$,
$$H_m(\BA{}\otimes_{A^e}(X\otimes_BY))=
Tor_m^{A^e}(A,A)=HH_m(A,A)\;.$$
Moreover,
$$H_m((Y\otimes_A\BA{}\otimes_AY)\otimes_{B^e}B)=
Tor_m^{B^e}(B,B)=HH_m(B,B)\;.$$
Therefore, also for two $k$-algebras $A$ and $B$ the {\em Hochschild homology}
is an invariant of the derived category.

\section{Symmetric algebras and the K\"ulshammer-$\kappa$}
\label{kappa}

Given a perfect field $k$
of characteristic $p$ and a symmetric $k$-algebra $A$
with symmetrizing bilinear form $(\;,\;)_A=(\;,\;)$,
K\"ulshammer defined in \cite{Kueldeutsch}
a mapping $\zeta_n=\zeta_n^A:Z(A)\lra Z(A)$ by
$(z,a^{p^n})_A=(\zeta_n(z),a)_A^{p^n}$,
for all $z\in HH^0(A)$ and $a\in HH_0(A)$.

The bilinear form
$(\;,\;)_A$ induces an identification $A\simeq Hom_k(A,k)$ as
$A\otimes_kA^{op}$-bimodules.

In \cite{rogquest} we showed that
$F^e(Hom_k(A,k))\simeq Hom_k(B,k)$ and therefore we get an induced
non degenerate symmetric bilinear form $(\;,\;)_B$ on $B$, making $B$
a symmetric algebra as well.
In \cite{kuelsquest} we showed that $\zeta_n$ on $A$ induces the
mapping $\zeta_n$ on $B$ by the derived equivalence $F$ and $F^e$,
if one chooses the symmetrizing form $(\;,\;)_B$ on $B$.

We would like to mention a consequence which is implied by
\cite{kuelsquest}.

\begin{Cor}\label{exponent}
Suppose $k$ is a perfect field of characteristic $p$,
suppose $G$ and $H$ are finite groups, $B_G$
is a block of $kG$ with defect group $D_G$ and suppose
$B_H$ is a block of $kH$ of defect group $D_H$.
If $D^b(kG)\simeq D^b(kH)$, then the exponent of $D_G$
and the exponent of $D_H$ coincide.
\end{Cor}

Proof. This follows from the fact that the ideals $im(\zeta_n^{B_G})$ and
$im(\zeta_n^{B_H})$ are mapped to each other by the isomorphism
$Z(B_G)\simeq Z(B_H)$, and that by a result due to
K\"ulshammer \cite[formulae (17), (47) and (78)]{Kuelprog},
the exponent of $D_G$ is the smallest integer $n$ so that
$im(\zeta_n^{B_G})=im(\zeta_{n+1}^{B_G})$, and likewise for $H$.
\hfill\dickebox

\begin{Rem}
Suppose $k$ is an algebraically closed field of characteristic $p$,
suppose $G$ and $H$ are finite groups, suppose
$R$ is a complete discrete valuation domain of
characteristic $0$ with $R/rad(R)=k$ and field of fractions $K$.
Suppose $B_G$ is a block of $RG$ with defect group $D_G$ and suppose
$B_H$ is a block of $RH$ of defect group $D_H$.
If $D^b(B_G)\simeq D^b(B_H)$, then the orders of the defect
groups coincide: $|D_G|=|D_H|$.
The argument is the following
construction\footnote{This argument arose during a discussion with
Gabriele Nebe in Aachen in October 2005.
I am very grateful for her kind hospitality and
for giving me this reference} of Cliff,
Plesken and Weiss \cite{Plesken}. Let $B$ be a block of $RG$ with defect
group $D$ and
let $\Lambda$ be the centre of $B$. Define inductively $\Lambda_0:=\Lambda$
and $$\Lambda_{i+1}:=\{x\in K\Lambda\;|\;x\cdot rad(\Lambda_i)\subseteq
rad(\Lambda_i)\}\;.$$
Cliff, Plesken and Weiss show for algebraically closed fields $k$
\cite[Theorem 3.4]{Plesken} that
$$min\{s\;|\;\Lambda_s=\Lambda_{s+1}\}=|D|.$$
\end{Rem}

K\"ulshammer gives a second, in some sense dual mapping $\kappa_n$
defined by the equation $(z^{p^n},a)=(z,\kappa_n(a))^{p^n}$,
for $z\in HH^0(A)$, and for $a\in HH_0(A)$.

\begin{Prop}\label{kappainvariant}
The mapping $\kappa_n^A:HH_0(A)\lra HH_0(A)$ is invariant under a
derived equivalence $F:D^b(A)\lra D^b(B)$ in the sense that under
the induced isomorphism $HH_0(F):HH_0(A)\lra HH_0(B)$
one has $$HH_0(F)\circ \kappa_n^A\circ HH_0(F)^{-1}=\kappa_n^B$$
if one chooses the induced bilinear form $(\;,\;)_B$ on $B$.
\end{Prop}

Proof.
Examining this relation,
$\kappa_n$ fits in the commutative diagram
$$
\begin{array}{cccccc}
A\otimes_{A\otimes_kA^{op}}A&\simeq&
((A\otimes_{A\otimes_kA^{op}}A)^*)^*&\simeq&
Hom_k(Hom_{D^b(A\otimes_kA^{op})}(A,A),k) \\
&&&&\downarrow
\mbox{\scriptsize $((Fr^k)_*)^n$}\\
\uar\kappa_n&&&&
Hom_k(Hom_{D^b(A\otimes_kA^{op})}(A,A),k^{(n)})\\
&&&&\uparrow
\mbox{\scriptsize $((\mu_p)^*)^n$}\\
A\otimes_{A\otimes_kA^{op}}A&\simeq&
((A\otimes_{A\otimes_kA^{op}}A)^*)^*&\simeq&
Hom_k(Hom_{D^b(A\otimes_kA^{op})}(A,A),k)
\end{array}
$$
where $Fr$ is the Frobenius automorphism on $k$, and where
$$\mu_p:Z(A)\ni z\mapsto z^p\in Z(A).$$
Using $Z(A)=Hom_{A^e}(A,A)$, the mapping $\mu_p$ corresponds to
$$Hom_{A^e}(A,A)\ni f\mapsto
\underbrace{f\circ f\circ\dots\circ f\circ f}_{\mbox{\scriptsize
$p$ factors}}
\in Hom_{A^e}(A,A)\;.$$
Since $F^e$ is a functor it is  clear that
$F^e(\mu_p^A)=\mu_p^B$.
The rest of the proof is exactly analogous to the one in \cite{kuelsquest}.

Hence $F^e$ maps $\kappa_n^A$ to $\kappa_n^B$ with respect to the
induced identification $B\simeq Hom_k(B,k)$ and the proof is
finished. \hfill\dickebox

\bigskip

This result will be generalised to higher Hochschild
homology in Theorem~\ref{derivedkappacup}. Though the proof there
covers the present proposition it seems to be useful to have a
short independent proof here.

We use K\"ulshammer's description of the image and the kernel of $\kappa_n$
to get a nice invariant of derived categories.
Let $P_n(ZA):=<z^{p^n}|\;z\in Z(A)>_{k-\mbox{\scriptsize space}}$.
Then
$$P_n(ZA)^\perp/KA=\{x\in A/KA|\;(z^{p^n},x)=0\;\forall z\in Z(A)\}$$
and
$$T_n(ZA)^\perp/KA=\{x\in A/KA|\;\forall\; z\in Z(A):
z^{p^n}=0\Rightarrow (z,x)=0\}.$$
Now, $T_n(ZA)^\perp/KA$ is a $Z(A)$-submodule of the $Z(A)$-module $A/KA$.
Indeed, for $x\in T_n(ZA)^\perp/KA$ and $y\in Z(A)$, one gets for any
$z\in Z(A)$ with $z^{p^n}=0$ also $(yz)^{p^n}=0$ and so, $(z,yx)=(yz,x)=0$ as
well.

\begin{Cor}
Let $F:D^b(A)\lra D^b(B)$ be an equivalence of standard type
between the derived categories of the symmetric $k$-algebras $A$ and $B$
over a perfect field $k$. Then
\begin{itemize}
\item
the isomorphism $HH_0(F):A/KA\lra B/KB$ maps
$P_n(ZA)^\perp/KA$ to $P_n(ZB)^\perp/KB$.
\item
the isomorphism $HH_0(F):A/KA\lra B/KB$ maps
$T_n(ZA)^\perp/KA$ to $T_n(ZB)^\perp/KB$ as submodules over the
centres of the algebras.
\end{itemize}
\end{Cor}

Proof. This is a consequence of Proposition~\ref{kappainvariant}, the
fact that our isomorphisms are functorial, and hence preserve the
natural structure of $A/KA=HH_0(A,A)$ as $HH^0(A,A)=Z(A)$-module,
and the fact that the first module $P_n(ZA)^\perp/KA$ is the
kernel of $\kappa_n^A$, whereas $T_n(ZA)^\perp/KA$ is the image of
$\kappa_n^A$, as was shown by
K\"ulshammer~\cite[(52),(53)]{Kuelprog}.\hfill\dickebox

\begin{Rem}
Since the analogous statement for the, in some sense dual, mapping $\zeta$
proved to be extremely useful for distinguishing derived categories,
such as the classification of tame domestic symmetric algebras by
Holm-Skowro\'nski \cite{Holmskow} or to fix some of the open parameters
in the derived equivalence classification of tame blocks of group rings
\cite{hz-tame}, there is quite some hope that this corollary
is as useful as was the method in \cite{kuelsquest}. Indeed, this new
quite sophisticated invariant is explicitly computable in case the
Hochschild homology is known as vector space, but additional structure
is lacking.
\end{Rem}

A first step for the generalisation is the following lemma.
Instead of taking degree $0$ homology, one takes higher
degree homology, and finds again a bilinear form.

\begin{Lemma}\label{pairing}
Let $A$ be a symmetric $k$-algebra with symmetrizing form $(\;,\;)$.
Then, there is a non degenerate bilinear form
$$(\;,\;)_m:HH^m(A,A)\otimes HH_m(A,A)\lra k$$
so that $(\;,\;)_0=(\;,\;)|_{Z(A)\otimes A/KA}$\;.
This form gives an isomorphism $$Hom_k(HH^m(A,A),k)\simeq HH_m(A,A)$$
and any such $k$-linear isomorphism induces a non degenerate bilinear
form.
\end{Lemma}

Proof. Since
$\left(-\otimes_{A\otimes A^{op}}A,Hom_k(A,-)\right)$ is an adjoint pair
of functors between $k$-vector spaces and $A\otimes A^{op}$-modules,
$$Hom_k(\BA\otimes_{A\otimes A^{op}}A,k)
\simeq Hom_{A\otimes A^{op}}(\BA,Hom_k(A,k))$$

Taking homology of these complexes, and using that
$Hom_k(-,k)$ is exact and contravariant, gives
\begin{eqnarray*}
Hom_k(HH_m(A,A),k)&\simeq&
Hom_k(H_m(\BA\otimes_{A\otimes A^{op}}A),k)\\
&\simeq&
H^m(Hom_k(\BA\otimes_{A\otimes A^{op}}A,k))\\
&\simeq&
H^m(Hom_{A\otimes A^{op}}(\BA,Hom_k(A,k)))\\
&=&HH^m(A,Hom_k(A,k))
\end{eqnarray*}

Since $A$ is symmetric, $A\simeq Hom_k(A,k)$ as
$A\otimes A^{op}$-modules and we get an isomorphism
$$HH^m(A,A)\stackrel{\varphi_n}{\lra} Hom_k(HH_m(A,A),k)$$
as $k$-vector spaces. Now, put for any $f\in HH^m(A,A)$ and
$x\in HH_m(A,A)$
$$(f,x)_m:=(\varphi_m(f))(x)\;.$$
It is clear that $(\;,\;)_m$ is bilinear since $\varphi_m$ is
$k$-linear. The form $(\;,\;)_m$ is non degenerate since $\varphi$
is an isomorphism. In case $m=0$ we find back the form $(\;,\;)$
which was used to identify $A$ with $Hom_k(A,k)$.

The last part of the statement is well known.

This finally proves the lemma.\hfill\dickebox

\section{The cup product $\kappa$}
\label{cupkappa}

Now, $HH^*(A,A)$ carries a natural graded commutative ring
structure given by cup product. We shall describe how this
allows to define, just
as in the construction of $\kappa_n$, a semilinear $p$-power map
$\mu^{(m)}_p$ and by this means a higher degree
$\kappa^{(m)}_n:HH_{p^nm}(A,A)\lra HH_{m}(A,A)$ with
$\kappa_n=\kappa_n^{(0)}$.

Let $m\in2\N$. Then, for any $n\in\N$ and any $x\in HH_{p^nm}(A,A)$
one gets a $k$-linear map
\begin{eqnarray*}
HH^{m}(A,A)&\lra& k\\
 f&\mapsto &Fr_p^{-n}\left((f^{p^n},x)_{p^nm}\right)\;.
\end{eqnarray*}
This is hence an element in $Hom_k(HH^m(A,A),k)$. Now, by
Lemma~\ref{pairing}
\begin{eqnarray*}
HH_m(A,A)&\simeq &Hom_k(HH^m(A,A),k)\\
x&\mapsto&\left(f\mapsto (f,x)_{m}\right)
\end{eqnarray*}
Hence,
for all ${n\in\N}$ and for all ${x\in HH_{p^nm}(A,A)}$ there is a unique
$\kappa_n^{(m)}(x)\in HH_m(A,A)$ so that for all $f\in HH^m(A,A)$
one has
$$\left(f^{p^n},x\right)_{p^nm}=
\left(\left(f,\kappa_n^{(m)}(x)\right)_m\;\right)^{p^n}.$$

We defined for all $n\in \N$ and all $m\in 2\N$
a mapping
$$\kappa_n^{(m),A}=\kappa_n^{(m)}:HH_{p^nm}(A,A)\lra HH_m(A,A)$$
so that $\kappa_n^{(0),A}=\kappa_n$. In case $A$ is clear from
the context, we denote $\kappa_n^{(m),A}=\kappa_n^{(m)}$.

\begin{Rem}\label{graded}
Observe that since the Hochschild cohomology ring is {\em graded
commutative}, for $p$ odd the mapping $$HH^{2m+1}(A,A)\ni f\mapsto
f^p\in HH^{p\cdot (2m+1)}(A,A)$$ is $0$ for $m\in \N$.
Of course, the $0$-mapping
is semilinear as well.
If $p=2$, then being graded commutative is just the same
as being commutative, and so the restriction on $m$ is not
necessary.
\end{Rem}

\begin{Theorem}\label{derivedkappacup}
Let $A$ be a finite dimensional symmetric $k$-algebra over the field $k$ of
characteristic $p>0$. Let $B$ be a second algebra so that
$D^b(A)\simeq D^b(B)$ as triangulated categories.
Let $p$ be a prime and let $m\in\N$.
Then, there is a standard equivalence $F:D^b(A)\simeq D^b(B)$,
and any such standard equivalence induces an isomorphism
$HH_m(F):HH_m(A,A)\lra HH_m(B,B)$ of all Hochschild homology
groups satisfying
$$HH_{m}(F)\circ\kappa_n^{(m),A}\circ HH_{p^nm}(F)^{-1}=
\kappa_n^{(m),B}\;.$$
\end{Theorem}

Proof. By Remark~\ref{graded} in case $m$ odd we may assume $p=2$.

Again let $\BA{}$ be the bar resolution. Then the complex computing the
Hochschild homology is $\BA{}\otimes_{A^e}A$. Since
$A$ is finite dimensional, we compute for any integer $\ell$
\begin{eqnarray*}
H_{\ell}(\BA{}\otimes_{A^e}A)
&\simeq& Hom_k(H^{\ell}(Hom_k(\BA{}\otimes_{A^e}A,k)),k)\\
&\simeq& Hom_k(H^{\ell}(Hom_{A^e}(\BA{}, Hom_k(A,k))),k)
\end{eqnarray*}
just using the standard adjointness formulas between hom and
tensor functors. Finally, using the bilinear form $(\;,\;)$ on $A$
we get
\begin{eqnarray*}
H_{\ell}(\BA{}\otimes_{A^e}A)
&\simeq& Hom_k(H^{\ell}(Hom_{A^e}(\BA{}, Hom_k(A,k))),k)\\
&\simeq& Hom_k(H^{\ell}(Hom_{A^e}(\BA{}, A)),k)\;.
\end{eqnarray*}
We discover
$$HH_{\ell}(A,A)\simeq Hom_k(Ext^\ell_{A^e}(A,Hom_k(A,k)),k)\simeq
Hom_k(HH^\ell(A,A),k)$$
and get diagram
$$
\begin{array}{ccccc}
HH_{p^nm}(A,A)&\stackrel{\simeq}{\lra}&
Hom_k(Ext^{p^nm}_{A^e}(A,Hom_k(A,k)),k)
&\stackrel{\simeq}{\lra}&Hom_k(HH^{p^nm}(A,A),k)\\
&&&&\mbox{\scriptsize $((\mu^{(m)}_p)^n)^*$}\dar
\phantom{\mbox{\scriptsize $((\mu^{(m)}_p)^n)^*$}}\\
\phantom{\kappa_n^{(m)}}\dar\kappa_n^{(m)}&&&&Hom_k(HH^{m}(A,A),k^{(n)})\\
&&&&\mbox{\scriptsize $(Fr_p^n)_*$}\uar
\phantom{\mbox{\scriptsize $Fr_p^n$}}\\
HH_{m}(A,A)&\stackrel{\simeq}{\lra}&
Hom_k(Ext^{m}_{A^e}(A,Hom_k(A,k)),k)
&\stackrel{\simeq}{\lra}&Hom_k(HH^{m}(A,A),k)
\end{array}
$$
which we easily see to be commutative by what we observed previously.

We need to show that applying $HH_*(F)$ (resp. $HH^*(F)$)
to the various mapping spaces for the Hochschild (co-)homology of $A$
gives the analogous mapping for $B$.

Let $X\in D^b(A\otimes B^{op})$ be a two-sided tilting complex
with inverse $Y\in D^b(B\otimes A^{op})$. We may and will assume that
$X$ and $Y$ are complexes being projective on the left and projective
on the right. Then, we may replace the left derived tensor product by
the ordinary tensor product. Let
$$F_X:=Y\otimes_A-\otimes_AX:
D^b(A\otimes A^{op})\lra D^b(B\otimes B^{op})$$
Then, $F_X$ is a triangle equivalence, and in particular, $F_X(\BA{})$
is a resolution of $B$ in $D^b(B\otimes B^{op})$. So, $F_X$
induces a commutative diagram
$$\begin{array}{ccc}
Hom_k(HH_{p^nm}(A,A),k)&\lla&Ext^{p^nm}_{A^e}(A,A)\\
\dar&&\dar\\
Hom_k(HH_{p^nm}(B,B),k)&\lla&Ext^{p^nm}_{B^e}(B,B)
\end{array}$$
where the bottom row is again given by the adjointness formula
between Hom and tensor functors.

Moreover, $F_X(Hom_k(A,k))= Hom_k(B,k)$ as was shown in \cite{rogquest}.
This implies that an isomorphism $A\simeq Hom_k(A,k)$
induces an isomorphism
$B\simeq Hom_k(B,k)$ as bimodules. Therefore, the induced diagram
$$\begin{array}{ccc}
Hom_k(Ext^{p^nm}_{A^e}(A,Hom_k(A,k)),k)&\lra &Hom_k(Ext^{p^nm}_{A^e}(A,A),k)\\
\dar F_X&&\dar F_X\\
Hom_k(Ext^{p^nm}_{B^e}(B,Hom_k(B,k)),k)&\lra &Hom_k(Ext^{p^nm}_{B^e}(B,B),k)
\end{array}$$
is commutative.

Now, we know that the cup product on $HH^*(A,A)$ is the composition of
mappings in $Ext^{*}_{A^e}(A,A)=Hom_{D^b(A\otimes A^{op})}(A,A[*])$.
Therefore, applying $F_X$ to any of this composition of mappings is
going to give again the composition of mappings on
$Ext^{*}_{B^e}(B,B)=Hom_{D^b(B\otimes B^{op})}(B,B[*])$.
Since the mapping $\mu_p^{(m)}$ is the $k$-linear dual of this, again
$F_X$ induces a commutative diagram
$$\begin{array}{ccc}
Hom_k(HH^{p^nm}(A,A),k)&\stackrel{\mu_p^{(m),A}}{\lra}
&Hom_k(HH^{m}(A,A),k^{(n)})\\
\dar F_X&&\dar F_X\\
Hom_k(HH^{p^nm}(B,B),k)&\stackrel{\mu_p^{(m),B}}{\lra}
&Hom_k(HH^{m}(B,B),k^{(n)})
\end{array}$$

Since $F_X$ acts on the contravariant variable of
$Hom_k(HH^{m}(A,A),k)$ and in the space of semilinear mappings
$Hom_k(HH^{m}(A,A),k^{(n)})$,
we get
\begin{eqnarray*}
Hom_k(F_X,k)\circ \lefteqn{Hom_k(HH^m(A,A),(\mbox{Fr}_k))=}\\
&=&
Hom_k(HH^m(A,A),(\mbox{Fr}_k))\circ Hom_k(F_X,k)\;.
\end{eqnarray*}
This shows the theorem.
\hfill\dickebox

\bigskip

We observe first properties analogous to those in
K\"ulshammer~\cite{Kuelprog}. For this put
$$T_n^{(m)}(HH^m(A,A)):=\{x\in HH^m(A,A)\;|\;x^{p^n}=0\}.$$

\begin{Prop}\label{properties}
Suppose $k$ is a perfect field of characteristic $p>0$
and that $A$ is a symmetric finite dimensional $k$-algebra. Then, denoting by
$\perp_m$ the orthogonality with respect to the pairing $(\;,\;)_m$,
\begin{enumerate}
\item $\kappa_n^{(m)}$ is $k$-semilinear,
\item $\kappa_{n+\ell}^{(m)}=\kappa_\ell^{(m)}\circ\kappa_{n}^{(p^\ell m)}$
\item $im(\kappa_n^{(m)})=\left(T_n^{(m)}\right)^{\perp_m}$.
\item $ker\;\kappa_n^{(m)}=\{x^{p^n}\;|\;x\in HH^m(A,A)\}^{\perp_{p^nm}}$
\end{enumerate}
\end{Prop}

Proof. The first statement comes from the construction in
the proof of Theorem~\ref{derivedkappacup} of $\kappa$ as
composition of semilinear mappings.

The second statement again is implied by the following argument.
\begin{eqnarray*}
(f,\kappa_{n+\ell}^{(m)}(x))_{m}&=&(f^{p^{n+\ell}},x)_{p^{n+\ell}m}\\
&=&(f^{p^{\ell}},\kappa_{n}^{(p^\ell m)}(x))_{p^{\ell}m}\\
&=&(f,\kappa_\ell^{(m)}(\kappa_{n}^{(p^\ell m)}(x)))_{m}
\end{eqnarray*}

The third statement is shown as follows:
The defining equation
$$(x^{p^n},y)_{p^nm}=\left((x,\kappa_n^{(m)}(y))_m\right)^{p^n}$$
and the fact that $(\;,\;)_{p^nm}$ is non degenerate
show that $im(\kappa_n^{(m)})^{\perp_m}=T_n^{(m)}$.
Since $(\;,\;)_m$ is non degenerate, we may take the orthogonal
spaces of these and get the result.

The fourth statement comes directly from the defining equation
$$(x^{p^n},y)_{p^nm}=\left((x,\kappa_n^{(m)}(y))_m\right)^{p^n}$$
as well.
\phantom{x}\hfill\dickebox

\begin{Rem}
We see that the kernel and the image of $\kappa_n^{(m)}$ are
very much linked to the set of
nilpotent elements of the Hochschild cohomology. Snashall and
Solberg conjectured \cite{SnashallSolberg} that
the Hochschild cohomology ring of any finite dimensional algebra
is finitely generated modulo the ideal generated by nilpotent elements.
\end{Rem}

\begin{Cor}
We get $dim(im(\kappa_n^{(m)}))=dim(HH^m(A,A))-dim(T_n^{(m)})$.
\end{Cor}

Proof. This is an immediate consequence of the third statement of
Proposition~\ref{properties}. \hfill\dickebox

\bigskip

In general even degree Hochschild cohomology rings of symmetric algebras
contain nilpotent elements, but are not necessarily entirely nilpotent.
As an example I refer to the article Erdmann and
Holm \cite[Section 4]{ErdHolm} where Hochschild cohomology
rings of self-injective Nakayama algebras, which includes the
Hochschild cohomology of Brauer tree algebras,
are computed. There nilpotent elements arise in even Hochschild
degrees, though the even degree Hochschild cohomology modulo the nilpotent
radical is not zero in general.
In particular, $\kappa_n^{(m)}$ is neither zero nor surjective in general.

\begin{Rem}
In the joint paper \cite{Nonsymmetric} with Bessenrodt and Holm we showed
that for the degree zero Hochschild homology one may pass from a
possibly non-symmetric algebra $A$ to its trivial extension $\T A$.
Rickard showed in \cite{Rickard}
that whenever the algebras $A$ and $B$ are derived
equivalent then also the trivial extension algebras $\T A$ and $\T B$
are derived equivalent.
In degree $0$ it is then possible to interpret the mappings $\kappa$ and
$\zeta$ on the degree $0$ (co-)homology of $\T A$ in terms of $A$ only.
One might ask if an analogous construction is possible for
$\kappa_n^{(m)}$ as well. The obvious fact that the Hochschild homology
of $A$ is a direct factor of the Hochschild homology of $\T A$ might give
a natural definition. Nevertheless, there are quite a number of
technical problems, such as the fact that Rickard gives a one-sided tilting
complex only whereas a two-sided complex is needed for our method.
Moreover, on a more practical level, in order to be
able to compute $\kappa_n^{(m)}$ via the trivial extension
method, one needs at the present stage at least parts of the multiplicative
structure of the Hochschild cohomology of $\T A$. Even for rather small
algebras $A$ its trivial extension $\T A$ usually will have
quite complicate cohomology.
A significant simplification is needed and at the moment I do not see
clearly how one can cope with these difficulties.
\end{Rem}

\section{Stasheff-Quillen's construction of the Gerstenhaber structure and the
Gerstenhaber $\kappa$}
\label{Gerstenhaber}

We recall first a most helpful construction of the Gerstenhaber bracket
appearing in a slightly implicit fashion in Quillen \cite{Quillen}
and very explicitly in Stasheff \cite{Stasheff}.
I learned the construction in discussions from Bernhard Keller
\cite[Section 4.7]{KellerGersten}. This construction shows that
the Gerstenhaber bracket can be defined using a homological
construction on the bar complex.
For the reader's convenience we give the the construction in some detail.

\subsection{Stasheff-Quillen's construction}
\label{kellerconstruction}

Let
$$Coder(\B(A),\B(A)):=
\{D\in End_{A\otimes A^{op}}(\B(A))| \Delta\circ D=
(id_{\B(A)}\otimes D+D\otimes id_{\B(A)})\circ \Delta\}$$
be the coderivations. Since $\B(A)$ is graded, $Coder(\B(A),\B(A))$ is graded
as well. Denote by $Coder^n(\B(A),\B(A))$  the degree $n$ coderivations.
The vector space $Coder(\B(A),\B(A))$ is a
graded Lie algebra with Lie bracket being the commutator.
Now (cf e.g. Stasheff \cite[Proposition]{Stasheff}),
$$Coder(\B(A),\B(A))\simeq Hom_{A\otimes A^{op}}(\BA,A)[1]\;.$$
The isomorphism is induced by composing an $f\in Coder(\B(A),\B(A))$
with the projection $\tau$ on the degree $1$ component $A$ of $\B(A)$
so that $\gamma_A(f):=\tau\circ f\in Hom_{A\otimes A^{op}}(\BA,A)[1]$.
So, there is a unique $d_A\in Coder^1(\B(A),\B(A))$ with $\tau\circ d_A=m_A$
for $m_A$ being the multiplication map $m_A:A\otimes_kA\lra A$.
$m_A$ being associative is equivalent to $d_A^2=0$.
Hence, $Coder(\B(A),\B(A))$ is a differential graded Lie algebra.

\begin{Lemma} \label{Gerstenhaberpower}
(Keller, personal communication)
\begin{itemize}
\item
Suppose $k$ is a field. Then,
$$D\in Coder^{2n+1}(\B(A),\B(A))\Rightarrow
D^{2}\in Coder^{2\cdot(2n+1)}(\B(A),\B(A)).$$
\item
Suppose $k$ is a field of characteristic $p>0$.
Then, $$D\in Coder^{2n}(\B(A),\B(A))\Rightarrow
D^{p}\in Coder^{2pn}(\B(A),\B(A)).$$
\end{itemize}
\end{Lemma}

Proof of the first statement:
\begin{eqnarray*}
\Delta\circ D^2&=&
(id_{\B(A)}\otimes D+D\otimes id_{\B(A)})^2\circ \Delta\\
&=&id_{\B(A)}\otimes D^2+(D\otimes id_{\B(A)})(id_{\B(A)}\otimes D)+
(id_{\B(A)}\otimes D)(D\otimes id_{\B(A)})+D^2\otimes id_{\B(A)}\\
&=&id_{\B(A)}\otimes D^2-D\otimes D+D\otimes D+D^2\otimes id_{\B(A)}
\end{eqnarray*}
Proof of the second statement:
\begin{eqnarray*}
\Delta\circ D^p&=&
(id_{\B(A)}\otimes D+D\otimes id_{\B(A)})^p\circ \Delta\\
&=&\left(id_{\B(A)}\otimes D^p+
\left(\sum_{j=1}^{p-1}{p\choose j}\cdot (D^j\otimes D^{p-j})\right)+
D^p\otimes id_{\B(A)}\right)\circ\Delta\\
&=&\left(id_{\B(A)}\otimes D^p+D^p\otimes id_{\B(A)}\right)\circ\Delta
\end{eqnarray*}
\hfill\dickebox

\bigskip

\begin{Lemma}\label{gehtnicht}
Let $k$ be a field of characteristic $p$. Let $D\in Coder^n(\B(A),\B(A))$.
\begin{enumerate}
\item If $p=2$ and $n\in\N$,
then the mapping $D\mapsto D^2$ induces a mapping
$$HH^{n+1}(A,A)\lra HH^{2n+1}(A,A)$$
\item
If $p>2$ and $n=2m\in 2\N$, then
the mapping $D\mapsto D^p$ induces a mapping
$$HH^{2m+1}(A,A)\lra HH^{2pm+1}(A,A)$$
\end{enumerate}
\end{Lemma}

Proof.
Let $D\in Coder^n(\B(A),\B(A))$. Then, $D^p\in Coder^{pn}(\B(A),\B(A))$.
The differential in the Hochschild cohomology complex $Hom_{A^e}(\BA,A)$
corresponds to the commutator $[d_A,-]$ where $Coder^1( \B(A),\B(A))\ni d_A$
comes from $\tau\circ d_A=m_A:A\otimes A\lra A$ being
the multiplication in the algebra $A$.
We would like to show that the $p$-power operation
induces a genuine operation on Hochschild cohomology.

For this, it is immediate that $$[d_A,D]=0\Rightarrow [d_A,D^p]=0.$$
Hence, the $p$-power operation induces an operation on the
cycles of the Hochschild cohomology complex.

We need to show moreover that for any $E\in Coder^{n-1}(A,A)$ one has
$$(D+[d_A,E])^p\in D^p+im([d_A,-]).$$
\begin{enumerate}
\item
Let $p=2$. Then, since $d_A^2=0$, we get $[d_A,E]^2=[d_A,E[d_A,E]]$
and therefore
$$(D+[d_A,E])^2=D^2+D[d_A,E]+[d_A,E]D+[d_A,E[d_A,E]].$$
We need to show that
$D[d_A,E]+[d_A,E]D\in im([d_A,-])$.
But, $$D[d_A,E]+[d_A,E]D=[d_A,[D,E]]+[[d_A,D],E]$$
and whenever $D$ is a Hochschild cocycle, then $[d_A,D]=0$ and therefore
$$(D+[d_A,E])^2=D^2+[d_A,[D,E]+E[d_A,E]]\;.$$
This shows the statement for $p=2$.
\item
Let $p>2$. Again using that $d_A^2=0$
one sees that for any positive integer $n$ one has
$$[d_A,E]^n=
(d_AE)^n-\left(\sum_{j=1}^{n-1}(d_AE)^j(Ed_A)^{n-j}\right)+(-1)^n(Ed_A)^n
$$
But, using as well that $d_A^2=0$,
for the $n$-fold Lie-bracket one gets the same result
$$
[d_A,E[d_A,\dots,E[d_A,E]\dots]]=
(d_AE)^n-\left(\sum_{j=1}^{n-1}(d_AE)^j(Ed_A)^{n-j}\right)+(-1)^n(Ed_A)^n
$$
So, there is an element $X_n(E)$ with $[d_A,E]^n=[d_A,X_n(E)]$ and hence
$[d_A,E]^p\in im([d_A,-]).$
Moreover, just as in the case $p=2$, for any Hochschild cocycle $D$
one has $[d_A,D]=0$ and one gets $$(D+[d_A,E])^p\in D^p+im([d_A,-])\;.$$
\end{enumerate}
This finishes the proof.\hfill\dickebox

\subsection{The $p$-restricted Lie structure and its derived invariance}

We recall the definition of a restricted Lie algebra
(cf e.g.\cite[Chapter 2 Section 1]{Strade}).

\begin{Def}
Let $L$ be a Lie algebra over a field $k$ of characteristic $p>0$.
Denote $$\left(ad\;(a\otimes X+b\otimes 1)\right)^{p-1}(a\otimes 1)=
\sum_{i=1}^{p-1}i\cdot s_i(a,b)\otimes X^i\in L\otimes k[X].$$
A mapping $[p]:L\lra L$ is called a {\em $p$-mapping} if
\begin{enumerate}
\item $ad\;a^{[p]}=(ad\;a)^p$
\item $(\alpha a)^{[p]}=\alpha^pa^{[p]}$
\item $(a+b)^{[p]}=a^{[p]}+b^{[p]}+\sum_{i=1}^{p-1}s_i(a,b)$

\end{enumerate}
A Lie algebra $L$ together with a $p$-mapping $[p]$ is then
called a {\em $p$-restricted Lie algebra}.
\end{Def}

The following proposition should be well known to the experts, but
I could not find a reference. So, I include the short proof.

\begin{Prop}
\begin{itemize}
\item
For any field $k$ of characteristic $2$ the differential graded Lie algebra
$Coder^{\ast}(\B(A),\B(A){})$ and its homology $HH^{\ast+1}(A,A)$
are $2$-restricted Lie algebras under the Gerstenhaber construction.
\item
For any field $k$ of characteristic $p>2$ the sum of the odd degree
Hochschild cohomology groups $\bigoplus_{n\in\N}HH^{2n+1}(A,A)$
is a $p$-restricted Lie algebra under the Gerstenhaber construction.
\end{itemize}
\end{Prop}

Proof.
As we have seen in Lemma~\ref{Gerstenhaberpower},
for $p=2$ the square of any coderivation
is a coderivation. Moreover by Lemma~\ref{Gerstenhaberpower},
for any prime $p$ the $p$-power of an even degree
coderivation is a coderivation. So, the Gerstenhaber $p$-power induces
a mapping $$Hom_{A^e}((\BA)^{2m+1},A)\lra Hom_{A^e}((\BA)^{2pm+1},A).$$

Let us show the first property $(ad\;a)^p=ad\;a^{[p]}$.
The Lie structure on the Hochschild cohomology is given by the commutator
bracket on the coderivations $Coder(\B(A),\B(A))$.
The mapping $a\mapsto a^{[p]}$ is given by taking
the ordinary composition of mappings.
Hence,
\begin{eqnarray*}
(ad\;a)^p(y)&=&
\underbrace{[a,[a,\dots,[a}_{\mbox{\scriptsize $p$ factors}},\;y
\underbrace{]\dots]]}_{\mbox{\scriptsize $p$ }}\\
&=&[a^{p},y]\\
&=&\left(ad\;a^{[p]}\right)(y)
\end{eqnarray*}

The second property is trivial, since
$(\alpha a)^p=\alpha^p a^p$ in $Coder(\B(A),\B(A))$.

The third property is done exactly analogously to
the first example following Lemma 1.2 in \cite[Chapter 2]{Strade}.
\hfill\dickebox

\begin{Prop}\label{Gerstenhaberkompatibel}
Let $A$ and $B$ be $k$-algebras over a field $k$.
Suppose $D^b(A)\simeq D^b(B)$
as triangulated categories.
\begin{itemize}
\item
If the characteristic of $k$ is $2$ then $HH^\ast(A,A)$ and
$HH^\ast(B,B)$ are isomorphic as restricted Lie algebras.
\item
If the characteristic of $k$ is $p>2$, then the Lie algebras
consisting of odd degree Hochschild cohomologies
$\bigoplus_{n\in\N}HH^{2n+1}(A,A)$ and $\bigoplus_{n\in\N}HH^{2n+1}(B,B)$
are isomorphic as restricted Lie algebras.
\end{itemize}
\end{Prop}

Proof. The fact that the Gerstenhaber structure is preserved is shown
by Keller in \cite{KellerGersten}.

We need to show that the isomorphism maps the $p$-power maps to each other.
Let $X\in D^b(A\otimes_kB^{op})$ be a twosided tilting complex
which we will assume to be formed by modules projective on either side.
Let $Y=Hom_k(X,k)\in D^b(B\otimes_kA^{op})$ be the inverse complex.

We first suppose $p=2$.
Then, as above, the functor
$$Y\otimes_A-\otimes_AX:D^b(A\otimes_kA^{op})\lra D^b(B\otimes_kB^{op})$$
is an equivalence. So, this functor induces an isomorphism
\begin{eqnarray*}
\varphi_X:Hom_{A^e}(\BA,A)&\lra&
Hom_{B^e}(Y\otimes_A\BA\otimes_AX,Y\otimes_AA\otimes_AX)\\
&&\;\;\;\;\;\;\|\\
&&Hom_{B^e}(\BB{},B)
\end{eqnarray*}
since $Y\otimes_A\BA\otimes_AX$ is a projective resolution of $B$ in
the category of $B\otimes B^{op}$-modules, and therefore homotopy equivalent
to $\BB{}$, and since $Y\otimes_AA\otimes_AX\simeq B$ in
$D^b(B\otimes_kB^{op})$.
We know by the discussion in Section~\ref{kellerconstruction}
that there is an isomorphism $\gamma_A$
\begin{eqnarray*}
Coder(\B(A),\B(A))[-1]&\stackrel{\gamma_A}{\lra}
&Hom_{A\otimes_kA^{op}}(\BA,A)\\
D&\mapsto&\tau\circ D
\end{eqnarray*}
where $\tau$ is the projection mapping $\B(A)\lra A$ on the degree $1$
component $A$ of $\B(A)$. Denote by $\sigma_A$ the square map on
$Coder(\B(A),\B(A))$ and likewise by $\sigma_B$ the squaring on
$Coder(\B(B){},\B(B){})$.

We need to show that
$$\gamma_B^{-1}\circ\sigma_B\circ\gamma_B=
\varphi_X\circ\gamma_A^{-1}\circ\sigma_A\circ
\gamma_A\circ\varphi_X^{-1}$$
or in other words that the diagram below is commutative.

{\small
\unitlength1cm
\begin{picture}(13,6)
\put(3,4){$Coder(\B(A){},\B(A){})[-1]$}
\put(8,4){$Coder(\B(A){},\B(A){})[-1]$}
\put(3,2){$Coder(\B(B){},\B(B){})[-1]$}
\put(8,2){$Coder(\B(B){},\B(B){})[-1]$}
\put(0,5){$Hom_{A^e}(\BA{},A)$}
\put(12,5){$Hom_{A^e}(\BA{},A)$}
\put(0,1){$Hom_{B^e}(\BB{},B)$}
\put(12,1){$Hom_{B^e}(\BB{},B)$}
\put(1,4.8){\vector(0,-1){3.4}}
\put(13,4.8){\vector(0,-1){3.4}}
\put(1.1,2.9){\scriptsize$\varphi_X$}
\put(13.1,2.9){\scriptsize$\varphi_X$}
\put(1.9,4.8){\vector(2,-1){1}}
\put(12,4.8){\vector(-2,-1){1}}
\put(1.9,1.4){\vector(2,1){1}}
\put(12,1.4){\vector(-2,1){1}}
\put(2.1,4.4){\scriptsize $\gamma_A$}
\put(11.8,4.4){\scriptsize $\gamma_A$}
\put(2.1,1.7){\scriptsize $\gamma_B$}
\put(11.7,1.7){\scriptsize $\gamma_B$}
\put(6.8,4.1){\vector(1,0){1.0}}
\put(7,4.2){\scriptsize $\sigma_A$}
\put(6.8,2.1){\vector(1,0){1.0}}
\put(7,2.2){\scriptsize $\sigma_B$}
\end{picture}}

But this is obvious.

Suppose now $p>2$. Then, the $p$-power operation $\sigma_p$
is only defined on the space
$Coder^{2\N}(\B(A){},\B(A){})$. Restricting therefore
to only odd degree Hochschild cocycles, the same proof then holds.

This finishes the proof of the Proposition.
\hfill\dickebox

\begin{Rem}\begin{enumerate}
\item
In \cite{Kellerhigher}
Bernhard Keller shows that a derived equivalence of standard type
even preserves the structure of $Hom_{A^e}(\BA{},A)$ as $B_{\infty}$-algebra.
\item
Since there is no obvious reason why the $p$-power mapping should be
additive in general, neither semilinear,
it seems to be difficult to get an analogue for Gerstenhaber structures
to K\"ulshammer's mappings~$\kappa_n$.
\end{enumerate}
\end{Rem}

\subsection*{Acknowledgement:}
I want to thank Bernhard Keller for his explanations concerning the
Gerstenhaber bracket and I want to thank Gabriele Nebe for her kind
hospitality at the RWTH Aachen as well as for bringing my attention
to \cite{Plesken}.

Furthermore I want to thank the referee for pointing out the references
\cite{Quillen}, \cite{Stasheff} and very helpful remarks.

\end{document}